\documentclass[10pt]{article}
\usepackage{amsmath,amssymb}

\newtheorem{theorem}{Theorem}[section]

\newtheorem{proposition}[theorem]{Proposition}

\newcommand{\qed}{\enspace\vrule  height6pt  width4pt  depth2pt}
\newenvironment{proof}{\par\noindent{\bf Proof.}}{$\qed$\par\bigskip}

\newcommand{\supp}{\operatorname{supp}}
\newcommand{\U}{\operatorname{U}}

\begin{document}
\date{}
\title{
Semigroup algebras of submonoids of
polycyclic-by-finite groups and maximal
orders\thanks{Research partially supported by the
Onderzoeksraad of Vrije Universiteit Brussel, Fonds
voor Wetenschappelijk Onderzoek (Flanders),
Flemish-Polish bilateral agreement BIL2005/VUB/06
and a MNiSW research grant N201 004 32/0088
(Poland).}}
\author{Isabel Goffa\footnote{Research funded by a Ph.D grant of
the Institute for the Promotion of Innovation
through Science and Technology in Flanders
(IWT-Vlaanderen).} \and Eric Jespers \and Jan Okni\'{n}ski}
\maketitle

\begin{center}
Dedicated to Fred Van Oystaeyen, on the occasion of
his sixtieth birthday
\end{center}

\begin{abstract}
Necessary and sufficient conditions are given for
a prime  Noetherian algebra $K[S]$ of a submonoid
$S$ of a polycyclic-by-finite group $G$ to be a
maximal order. These conditions are entirely in
terms of the monoid $S$. This extends  earlier
results of Brown concerned with the group ring
case and of the authors for the case where $K[S]$
satisfies a polynomial identity.
\end{abstract}

\section{Introduction}

In this paper we continue our  investigations in
\cite{gof-jes,PI,okni}  on semigroup algebras $K[S]$
that are prime Noetherian maximal orders (for a
survey we refer the reader to \cite{boekjesok}).
There are two main issues to be dealt with. First,
when such algebras are Noetherian and second when
they are a maximal order. We briefly give some
background.  Recall that group algebras of
polycyclic-by-finite groups are the only known
examples of Noetherian group algebras. In
\cite{brown11,brown12}, K.A. Brown characterized
when such group algebras are prime maximal orders.
In the search for more classes of prime Noetherian
maximal orders, it is thus natural to consider
subalgebras of Noetherian group algebras. In
\cite{jesokn-noe} it is proved that the semigroup
algebra $K[S]$ of a submonoid $S$ of a
polycyclic-by-finite group is right Noetherian if
and only if $S$ has a group  of right quotients
$G=SS^{-1}$, with normal subgroups $F$ and $N$  such
that $F\subseteq S\cap N$, $G/N$ is finite, $N/F$ is
abelian and $S\cap N$ is finitely generated. In
particular, in this situation,  $S$ is finitely
generated and, if the unit group $\U(S)$ is trivial
then $K[S]$ satisfies a polynomial identity.
Furthermore, it follows that such a semigroup
algebra is right Noetherian if and only if it is
left Noetherian. We simply call such algebras
Noetherian. In \cite{PI}, the authors determined
conditions under which the semigroup algebra of a
submonoid $S$ of a finitely generated
abelian-by-finite group is a prime Noetherian
maximal order. It turns out  that the action of the
group of quotients on the minimal primes of some
abelian submonoid of $S$ is very important (some
invariance condition is crucial).

In this paper we deal with the general case,
provided such an invariance condition holds. Crucial
for our investigations is Theorem~1.1 in \cite{PI}
that says that the height one prime ideals $P$  of a
prime Noetherian algebra $K[S]$ with $P\cap S\neq
\emptyset$, where $S$ is a submonoid of a
polycyclic-by-finite group, are precisely the ideals
of the form $K[Q]$ with $Q$ a minimal prime ideal of
$S$ (recall that the other height one prime ideals
are contractions of height one prime ideals of the
group algebra $K[SS^{-1}]$). The set of all minimal
primes of $S$ will be denoted by $X^{1}(S)$.

For an arbitrary abelian monoid $A$, Anderson
\cite{and1,and} (see also \cite{cho,gil}) proved
that $K[A]$ is a prime Noetherian maximal order
if and only if $A$ is a finitely generated
submonoid of a torsion free  abelian group, so
that $A$ is a maximal order in its group of
quotients. In that case, these monoids $A$ are
precisely the finitely generated abelian monoids
$A$ so that $A=\U(A)\times A_{1}$ where
$A_{1}=A_{1}A_{1}^{-1}\cap F^{+}$ with $F^{+}$ a
positive cone of a free abelian group $F$ that
contains the group of quotients $A_{1}A_{1}^{-1}$
of a submonoid $A_{1}$ of $A$.

\section{Algebras of submonoids of polycyclic-by-finite groups}

Let $S$  be a submonoid of a polycyclic-by-finite
group  such that the semigroup algebra $K[S]$ is
Noetherian. Hence, $S$ has a group of quotients
$G=SS^{-1}$ with normal subgroups $F$ and $N$ such
that $F\subseteq S\cap N$, $N/F$ is abelian, $G/N$
is finite and $S\cap N$ is finitely generated.
Without loss of generality we may assume that the
groups $N$ and $N/F$ are torsion free. By $\sim_{F}$
we denote the congruence on $S$ defined by: $s
\sim_{F} t$ if and only if $s = ft$ for some $f \in
F$. The set of $\rho$-classes $S/ \rho$ has a
natural semigroup structure inherited from $S$ and
we denote this by $S/F$. Because of the natural
bijection between the minimal primes of $S$ and the
minimal primes of $S/F$, it is easily verified that
$S$ is a maximal order in its group of quotients $G$
if and only if the semigroup $S/F$ is a maximal
order in its group of quotients $G/F$. Throughout
this paper we will freely use all this notation.
Recall from \cite[Lemma 4.1.3]{boekjesok} that also
$K[S\cap N]$ is Noetherian. Moreover $K[S\cap N]$ is
a domain as $N$ is torsion free (see \cite[Theorem
37.5]{passman}).

 The following notation will be used.  For an element
$\alpha =\sum_{s\in S} k_{s}s\in K[S]$, with each
$k_{s}\in K$, we put $\supp (\alpha )=\{ s\in S\mid
k_{s}\neq 0\}$, the support of $\alpha$. By
$Q_{cl}(R)$ we denote the classical ring of
quotients of a prime Noetherian ring $R$. Recall
that $R$ is said to be a maximal order if the
following property holds  for every subring $T$ of
$Q_{cl}(R)$ with $R\subseteq T$: if there exist
regular elements $r_{1},r_{2}\in R$ so that
$r_{1}Tr_{2}\subseteq R$ then $R=T$. Equivalently,
$(I:_{r}I)=(I:_{l}I)=R$ for every (fractional) ideal
$I$ of $R$; here we put $(I:_{r}I)=\{ q\in Q_{cl}(R)
\mid Iq\subseteq I\}$ and similarly one defines
$(I:_{l}I)$. For more information and details we
refer the reader to \cite[Section 3.6]{boekjesok}.
In order to prove the main theorem we need the
following proposition.

\begin{proposition}\label{0.2}
Let $S$ be a submonoid of a polycyclic-by-finite
group $G$ such that the semigroup algebra $K[S]$ is
Noetherian. Then the semigroup $S\cap N$ is a
maximal order in its group of quotients if and only
if the semigroup algebra $K[S\cap N]$ is a maximal
order in its classical ring of quotients (with $N$ a
torsion free subgroup of finite index in
$G=SS^{-1}$, as above).
\end{proposition}
\begin{proof}
If $K[S\cap N]$ is a prime maximal order, then it is
well known and easy to prove that the semigroup $S
\cap N$ is a maximal order. Conversely, assume
$S\cap N$ is a maximal order in its group of
quotients $N$. Because $N$ is torsion free, we know
that $K[S\cap N]$ is a Noetherian domain. To prove
that $K[S\cap N]$ is a maximal order, let $I$ be a
non-zero ideal of $K[S\cap N]$ and let $0\neq q \in
Q_{cl}(K[S\cap N])$ be such that $qI \subseteq I$.
Then $q IK[N] \subseteq IK[N]$. Note that $IK[N]$ is
a two-sided ideal of $K[N]$ (see for example
\cite[Theorem~9.20]{goodenwar}). Because $N$ is a
torsion free polycyclic-by-finite group, we know
from Brown's result \cite{brown11} that $K[N]$ is a
maximal order. Hence, it follows that $q \in K[N]$.
Now, since $K[N] = K[F] \ast (N/F)$, a crossed
product of the finitely generated torsion free
abelian group $N/F$ over the group algebra $K[F]$,
we have that $q = \sum_{i = 1}^{n} \alpha_{i}
q_{i}$, with $\alpha_{i} \in K[F]$ and all $q_{i}$
are in  a transversal of $F$ in $N$. The image of
$q_{i}$ in $N/F$ we denote by $\overline{q_{i}}$.
Let $\leq $ denote an ordering on the group $N/F$.
Then, we may assume that $\overline{q_{1}} < \cdots
< \overline{q_{n}}$. Every nonzero element of $I$
can be written in the form $\beta t + \alpha $ for
some $t\in S\cap N, 0\neq \beta \in K[F]$ and some
$\alpha\in K[S\cap N]$ such that $\overline{s}<
\overline{t}$ for all $s\in \supp (\alpha)$ (if
$\alpha \neq 0$). Let $h(I)$ denote the set
consisting of all such possible elements $t\in S\cap
N$. Then $h(I)$ is an ideal of $S \cap N$. Since $q
I \subseteq I$, $N/F$ is ordered and $K[N]$ is a
domain, we get (using a standard graded algebra
argument) that $q_{n}(h(I)) \subseteq h(I)$. As
$S\cap N$ is a maximal order, this implies that
$q_{n} \in S\cap N$ and thus $\alpha_{n} q_{n} \in
(I:_{l}I)\cap K[S\cap N]$. So, $q-\alpha_{n}q_{n}\in
(I:_{l}I)$ and $|\supp (q-\alpha_{n}q_{n})|< |\supp
(q)|$. Hence, by an induction argument, we may
assume that $q-\alpha_{n}q_{n}\in K[S\cap N]$. So
$q\in K[S\cap N]$, as desired.

Similarly, one shows that $(I :_{r} I) = K[S\cap
N]$.
\end{proof}

In order to state the  main result we need some more
notation. By $\Delta^{+}(G)$ we denote the torsion
subgroup of the finite conjugacy center $\Delta(G)$
of a group $G$. It is well known that $K[G]$ is
prime if and only if $\Delta^{+}(G) = \{1\}$ (see
\cite[Theorem~5.5]{pas}). The following terminology
is used in \cite{brown11}. A group $G$ is said to be
dihedral-free if the normalizer of any subgroup $H$
isomorphic with the infinite dihedral group is of
infinite index in $G$.

\begin{theorem}\label{firstp-b-f}
Let $S$ be a submonoid of a polycyclic-by-finite
group such that the semigroup algebra $K[S]$ is
Noetherian, i.e., there exist normal subgroups $F$
and $N$ of $G = SS^{-1}$ such that $F\subseteq S\cap
N$, $N/F$  is abelian, $G/N$ is finite and $S\cap N$
is finitely generated. Suppose that for every
minimal prime $P$ of $S$ the intersection $P \cap N$
is $G$-invariant.

Then, the semigroup algebra $K[S]$ is a  prime
maximal order if and only if the monoid $S$ is a
maximal order in its group of quotients $G$, the
group $G$ is dihedral-free and $\Delta^{+}(G) =
\{1\}$.
\end{theorem}

\begin{proof}
First note that the $G$-invariance of $P\cap N$, for
every minimal prime ideal $P$ of $S$, is inherited
on $P\cap M=(P\cap N)\cap M$, for any normal
subgroup $M$ of $G$ with $M\subseteq N$ and $N/M$
finite. In particular, we may for the remainder
assume that $N$ is torsion free (and $N/F$ is
torsion free).

 If $K[S]$ is a prime maximal order
then (as before) $S$ is a maximal order in $G$.
Furthermore, because $K[G]$ is a localization of
$K[S]$, we know that $K[G]$ is a prime maximal order
as well. Hence, by Brown's result \cite{brown11}, it
follows that the the group $G$ is dihedral-free and
$\Delta^{+}(G) = \{1\}$. For the converse
implication, suppose that $S$ is a maximal order in
its group of quotients $G$, the group $G$ is
dihedral-free and $\Delta^{+}(G) = \{1\}$. Hence,
$K[G]$ and therefore also $K[S]$ is prime. Again, by
Brown's result, $K[G]$ is a maximal order.

Let $P$ be a minimal prime ideal of $S$. Then, by
\cite[Theorem 1.1]{PI} (see the introduction),
$K[P]$ is a height one prime of $K[S]$. Clearly,
$K[S]$ has a natural $G/N$-gradation with
homogeneous component of degree $e$ (the identity of
$G/N$) the algebra $K[S\cap N]$. So,  from
\cite[Theorem~17.9]{passman}, it then follows that
$$\widetilde{P}(N) = K[P]\cap K[S\cap N] = K[P\cap
N]= Q_{1}\cap \cdots \cap Q_{n},$$ with each $Q_{i}$
a height one prime ideal of $K[S\cap N]$; and these
are all the height one primes of $K[S\cap N]$
containing $\widetilde{P}(N)$. Because of the
assumption on the invariance of $P\cap N$, it easily
is verified that the set $\{Q_{1},\ldots ,Q_{n}\}$
is a full orbit (under the conjugation action) of
height one primes in $K[S\cap N]$. Clearly
$Q_{i}\cap (S\cap N)\neq \emptyset$. So, again by
the result mentioned in the introduction, $Q_{i} =
K[Q_{i}\cap (S \cap N)]$; moreover, $Q_{i}\cap
(S\cap N)$ is a minimal prime ideal of $S\cap N$ and
these are all the minimal primes of $S\cap N$
containing $P\cap N$. Because $K[S\cap N]$ is
Noetherian and $N$ is a polycyclic-by-finite group,
we also know (see for example \cite[Corollary
4.4.12]{boekjesok}) that each $Q_{i}$ contains a
normal element $n_{i}$, that is an element such that
$(S\cap N)n_{i}=n_{i}(S\cap N)$. Furthermore,
because $S$ is a maximal order, we get that   $S/F$
is a maximal order in its group of quotients $G/F$,
in which $N/F$ is abelian and of finite index.
Hence, it follows from \cite[Lemma~7.1.1]{boekjesok}
that $(S\cap N)/F$ is a maximal order as well.
Consequently, $S\cap N$ is a maximal order.  So, by
Proposition~\ref{0.2}, the Noetherian algebra
$K[S\cap N]$ is a maximal order. Since each $Q_{i}$
is a height one prime containing a divisorial ideal
(namely $K[(S\cap N)n_{i}]$), it therefore follows
that it is a divisorial height one prime ideal. It
then follows from Proposition~1.9 and
Proposition~1.10 in \cite{chamarie-a} that each
$Q_{i}$ is localizable (the localization will be
denoted  $K[S\cap N]_{Q_{i}}$) and hence also that
$\widetilde{P}(N)$ is a localizable semiprime ideal
of $K[S\cap N]$. Furthermore, $K[S\cap
N]_{\widetilde{P}(N)} = K[S\cap N]_{Q_{1}}\cap
\cdots \cap K[S\cap N]_{Q_{n}}$. Here we denote by
$K[S\cap N]_{\widetilde{P}(N)}$ the localization of
$K[S\cap N]$ with respect to the set $C_{N}(P) = \{c
\in K[S\cap N]\mid c + K[P\cap N] \mbox{ is a
regular element of the ring } K[S\cap N] / K[P\cap
N] \}$. Moreover (see for example
\cite[Lemma~13.3.5]{pas}), $C_{N}(P)$ is an Ore set
of regular elements of $K[S]$ and thus an element
$c\in K[S\cap N]$ belongs to $C_{N}(P)$ if and only
if $c+K[P]$ is regular in $K[S]/K[P]$. We begin by
showing that the localized ring $K[S]_{C_{N}(P)}$ is
a maximal order. To do so, we show that
$K[S]_{C_{N}(P)}$ is a local ring with unique
maximal ideal $P K[S]_{C_{N}(P)}$ and so that $P
K[S]_{C_{N}(P)}$ is invertible and every proper
non-zero ideal of $K[S]_{C_{N}(P)}$ is of the form
$(P K[S]_{C_{N}(P)})^{k}$ for some positive integer
$k$.

From \cite[Th\'eor\`eme~4.1.6]{cham} we know that
$K[S\cap N]_{\widetilde{P}(N)}$ is a semi-local
maximal order that is a principal left and right
ideal ring. In particular, it is an Asano order, it
has dimension one and its Jacobson radical is equal
to $\widetilde{P}(N) K[S\cap N]_{\widetilde{P}(N)}$.
Since this ring is the component of degree $e$ of
the $G/N$-graded ring $K[S]_{C_{N}(P)}$, it follows
from \cite[Theorem~17.9]{passman} that
$K[S]_{C_{N}(P)}$ also has dimension one. By the
above, it follows that the non-zero prime ideals of
$K[S\cap N]_{\widetilde{P}(N)}$ are precisely the
$n$ prime ideals $Q_{i}K[S\cap N]_{\widetilde{P}(N)}
= K[Q_{i}\cap S\cap N] K[S\cap
N]_{\widetilde{P}(N)}$. Consequently, the height one
primes of $K[S\cap N]_{\widetilde{P}(N)}$ are the
form $K[S\cap N]_{\widetilde{P}(N)} I$, where $I$ is
an ideal of $S\cap N$. Since $K[S\cap
N]_{\widetilde{P}(N)}$ is an Asano order, all its
non-zero ideals are products of height one prime
ideals. Hence all non-zero ideals of $K[S\cap
N]_{\widetilde{P}(N)}$ are of the form $K[S\cap
N]_{\widetilde{P}(N)} I$, where $I$ is an ideal of
$S\cap N$.  We now show that $PK[S]_{C_{N}(P)}$ is
the only height one prime of $K[S]_{C_{N}(P)}$. So,
let $Q$ be a height one prime ideal of
$K[S]_{C_{N}(P)}$. Again, because of the
$G/N$-gradation, $I=Q\cap K[S\cap
N]_{\widetilde{P}(N)}=K[S\cap N\cap
Q]_{\widetilde{P}(N)}$ is a non-zero semiprime ideal
of $K[S\cap N]_{\widetilde{P}(N)}$. Hence,
as explained above, non-zero ideals of $K[S\cap
N]_{\widetilde{P}(N)}$ are generated by their
intersection with $S\cap N$, and $(Q\cap S)\cap N =I\cap
(S\cap N)$ is an intersection of some of the
$Q_{i}\cap (S\cap N)$. In particular, $Q\cap S\neq
\emptyset$. Clearly, $Q\cap K[S]$ is a height one
prime ideal of $K[S]$ and thus $Q\cap K[S]=K[Q\cap
S]$ and $Q\cap S$ is a minimal prime ideal of $S$.

The assumptions therefore imply that $Q\cap (S\cap
N)$ is $G$-invariant. Since $Q\cap (S\cap
N)\subseteq Q_{i}$, for some $i$, we thus obtain
that $$Q\cap(S\cap N)\subseteq
\bigcap_{i=1}^{n}Q_{i} = K[P\cap N].$$ From
\cite[Proposition~1.3]{PI}, it follows that $$P =
{\mathcal B}(S((Q_{1}\cap (S\cap N))\cap \cdots \cap
(Q_{n}\cap (S\cap N)))S) = Q\cap S,$$  where
${\mathcal B}(J)$ denotes the prime radical of an
ideal $J$ of $S$. It follows that $Q =
PK[S]_{C_{N}(P)}$, as desired.

As, by assumption, $S$ is a maximal order, it
follows that $P(S : P)$ is an ideal of $S$ that is
not contained in $P$. Hence $P K[S]_{C_{N}(P)} (S :
P)$ is an ideal of $K[S]_{C_{N}(P)}$ that is not
contained in $P K[S]_{C_{N}(P)}$. Consequently, $P
K[S]_{C_{N}(P)}(S : P) = K[S]_{C_{N}(P)}$, i.e.
$PK[S]_{C_{N}(P)}$ is invertible. In particular, by
\cite[Proposition~4.2.6]{conandrob} this ideal
satisfies the Artin-Rees property. So, by a result
of P. Smith (see \cite[Theorem 11.2.13]{pas}),
$\bigcap_{k}(PK[S]_{C_{N}(P)})^{k} = \{0\}$. It then
easily follows (and it is well known) that every
proper non-zero ideal of $K[S]_{C_{N}(P)}$ is of the
form $(P K[S]_{C_{N}(P)})^{k}$, for some unique
positive integer $k$. So, each non-zero ideal of
$K[S]_{C_{N}(P)}$ is invertible. This proves the
desired properties of $K[S]_{C_{N}(P)}$ and thus
$K[S]_{C_{N}(P)}$ is a maximal order.

Next we will prove that  $\bigcap_{P\in X^{1}(S)}
K[S]_{C_{N}(P)}\cap G=S$. For this, suppose $g \in
\bigcap_{P\in X^{1}(S)} K[S]_{C_{N}(P)}\cap G$.
Then, for every minimal prime $P$ of $S$, there
exists an element $\beta_{P} \in C_{N}(P)$ such that
$\beta_{P}g \in K[S]$. We can assume that the image
$\overline{\beta_{P}}$ of $\beta_{P}$ is central in
$K[G/F]$. Indeed, since $\widetilde{P}(N)$ is
$G$-invariant, it follows that $\prod_{g\in
T}\beta_{P}^{g} \in C_{N}(P)$ (product in any fixed
order) for some finite transversal $T$ for $N$ in
$G$. Since $N/F$ is abelian, it follows that
$\overline{\prod_{g\in T}\beta_{P}^{g}}$ is central
in $K[G/F]$ and we can replace
$\overline{\beta_{P}}$ by this product. Furthermore,
$\overline{\beta_{P}}\overline{S}\overline{g}
\subseteq K[S/F] = K[\overline{S}]$ and hence
$\overline{S}\supp(\overline{\beta_{P}})\overline{S}\overline{g}
\subseteq \overline{S}$. The union
$\bigcup_{P}\overline{S}
\supp(\overline{\beta_{P}})\overline{S}$, over all
the minimal primes $P$ of $S$, is an ideal of
$\overline{S}$ that is not contained in any minimal
prime $\overline{P}$ of $\overline{S}$. Since also
$\bigcup_{P}
\overline{S}\supp(\overline{\beta_{P}})\overline{S}
\overline{g} \subseteq \overline{S}$, and, since
$\overline{S}$ is a maximal order by the comment in
the beginning of this section, it follows that
$\overline{g} \in \overline{S}$.  Hence $g\in S$, as
desired.

Because of the remark stated in the beginning of the
proof, the previous holds for any normal subgroup
$M$ of $G$ with $M\subseteq N$ and $N/M$ finite.

Finally, we prove the following claim: $K[S] =
\bigcap_{P\in X^{1}(S),M} K[S]_{C_{M}(P)}\cap K[G]$,
with $M$ running through all  torsion free normal
subgroups of $G$ with $M\subseteq N$ and $N/M$
finite. Note that this claim implies the result,
i.e., $K[S]$ is a maximal order. Indeed, let $I$ be
a non-zero ideal of $K[S]$ and suppose that $q\in
Q_{cl}(K[S])$ is such that $q I \subseteq I$. Let
$P$ be a minimal prime of $S$ and let $M$ be a
subgroup as described. Since $q I \subseteq I$, we
get $q I K[S]_{C_{M}(P)} \subseteq I
K[S]_{C_{M}(P)}$, with $I K[S]_{C_{M}(P)}$ a
two-sided ideal of $K[S]_{C_{M}(P)}$ by
\cite[Theorem~9.20]{goodenwar}. As $K[S]_{C_{M}(P)}$
is a maximal order, this yields $q \in
K[S]_{C_{M}(P)}$. On the other hand, as also $qIK[G]
\subseteq IK[G]$ and as $K[G]$ is a maximal order,
we get that $q \in K[G]$. Hence the claim implies
that $q \in K[S]$. So we have shown that
$(I:_{l}I)=K[S]$. Similarly, $(I:_{r}I)=K[S]$, and
thus indeed $K[S]$ is a maximal order.

So, to prove the claim,  let $q = \sum_{i= 1}^{n}
k_{i}g_{i} \in \bigcap_{P\in X^{1}(S),M}
K[S]_{C_{M}(P)}\cap K[G]$, where $k_{i}\neq 0 \in K$
and $g_{i}\in G$ for each $1 \leq i \leq n$ and
$g_{i}\neq g_{j}$ for $i \neq j$. It is enough to
show that $q \in K[S]$. We prove this by induction
on $n$. If $n=1$ then $q = kg$ with $g \in
\bigcap_{P\in X^{1}(S),M} K[S]_{C_{M}(P)}\cap G$ and
it follows from the above that $g\in S$, as desired.
Hence assume $n
> 1$.

Because $G$ is residually finite, there exists a
normal subgroup of finite index $M_{0}$ in $G$ such
that $M_{0}\subseteq N$ and $g_{i} g_{j}^{-1} \not
\in M_{0}$ for all $i \neq j$. Note that
$C_{M_{1}\cap M_{2} }(P) \subseteq C_{M_{1}}(P)
\subseteq C_{N}(P)$ for any two normal subgroups
$M_{1},M_{2}$ of $G$ so that $M_{1},M_{2}\subseteq
N$ and each $N/M_{i}$ is finite. Hence, in the
intersection $\bigcap_{P\in X^{1}(S),M}
K[S]_{C_{M}(P)}\cap K[G]$ we may assume that $M$
runs through all normal subgroups of $G$ with
$M\subseteq M_{0}$ and $M_{0}/M$ finite. In other
words we may replace $N$ by $M_{0}$ in the
intersection. It follows that the intersection
$\bigcap_{P\in X^{1}(S),M} K[S]_{C_{M}(P)}\cap K[G]$
is a $G/M_{0}$-graded ring. Hence, the induction
hypothesis yields that we may assume that $q$ is
$G/M_{0}$-homogeneous, that is, each
$g_{i}g_{j}^{-1}\in M_{0}$. Consequently, $n = 1$
and thus by the above we get $q \in K[S]$. This ends
the proof.
\end{proof}

Suppose that in the  previous  theorem one also
assumes that the group $G=SS^{-1}$ is
abelian-by-finite. Then, in \cite{PI}, it is
shown that  the condition ``for every minimal
prime $P$ of $S$ the intersection $P \cap N$ is
$G$-invariant'' is necessary for $K[S]$ to be a
maximal order. It is unknown whether this
necessity holds in general, nor it is known
whether this condition is redundant. That is, no
example of a  maximal order $S$ in a
polycyclic-by-finite group $G$ (with
$\Delta^{+}(G)=\{ 1\}$ and $G$ dihedral-free) is
known so that $K[S]$ is not a maximal order.
Proposition~\ref{0.2} shows that if such a monoid
$S$ exists then $G$ does not contain a normal
subgroup $F$ so that $F\subseteq \U(S)$ and $G/F$
is torsion free abelian.

\vspace{20pt}

 \noindent
 \begin{tabular}{ll}
 I. Goffa and E. Jespers & J. Okni\'{n}ski\\
 Department of Mathematics& Institute of Mathematics\\
 Vrije Universiteit Brussel & Warsaw University\\
 Pleinlaan 2& Banacha 2\\
 1050 Brussel, Belgium& 02-097 Warsaw, Poland\\
 efjesper@vub.ac.be and igoffa@vub.ac.be & okninski@mimuw.edu.pl
 \end{tabular}

\end{document}